# Optimally Tuned Iterative Reconstruction Algorithms for Compressed Sensing

Arian Maleki
Department of Electrical Engineering
Stanford University
arianm@stanford.edu

David L. Donoho
Department of Statistics
Stanford University
donoho@stanford.edu

*Abstract*— We conducted an extensive computational experiment, lasting multiple CPU-years, to optimally select parameters for two important classes of algorithms for finding sparse solutions of underdetermined systems of linear equations. We make the optimally tuned implementations available at `sparselab.stanford.edu`; they run 'out of the box' with no user tuning: it is not necessary to select thresholds or know the likely degree of sparsity.

Our class of algorithms includes iterative hard and soft thresholding with or without relaxation, as well as CoSaMP, subspace pursuit and some natural extensions. As a result, our optimally tuned algorithms dominate such proposals.

Our notion of optimality is defined in terms of phase transitions, i.e. we maximize the number of nonzeros at which the algorithm can successfully operate. We show that the phase transition is a well-defined quantity with our suite of random underdetermined linear systems. Our tuning gives the highest transition possible within each class of algorithms. We verify by extensive computation the robustness of our recommendations to the amplitude distribution of the nonzero coefficients as well as the matrix ensemble defining the underdetermined system.

Our findings include: (a) For all algorithms, the worst amplitude distribution for nonzeros is generally the constant-amplitude random-sign distribution, where all nonzeros are the same amplitude. (b) Various random matrix ensembles give the same phase transitions; random partial isometries may give different transitions and require different tuning; (c) Optimally tuned subspace pursuit dominates optimally tuned CoSaMP, particularly so when the system is almost square.

## I. INTRODUCTION

A recent flood of publications offers numerous schemes for obtaining sparse solutions of underdetermined systems of linear equations; a long list of useful ideas and suggestions can be gleaned from the papers [1]-[29], with new proposals appearing regularly. Popular methods have been developed from many viewpoints: $\ell_1$-minimization [1], [2], [3], [4], [5], [6], matching pursuit [7], [8], [9], [10], iterative thresholding methods [11], [12], [13], [14], [15], [16], [17], [18], [19], [20], subspace methods [10], [21], [20], [22], convex regularization [23], [24] and nonconvex optimization [25], [26], [27]. The specific proposals are often tailored to different viewpoints, ranging from formal analysis of algorithmic properties [12], [9], [28], [10], [6], to particular application requirements [13], [14], [15]. Such algorithms have potential applications in fields ranging from medical imaging to astronomy [29], [30].

The potential user now has a bewildering variety of ideas and suggestions that might be helpful, but this, paradoxically, creates uncertainty and may cause said potential user of such algorithms to avoid the topic entirely. In this note we announce a solution to this problem, in the form of freely available, optimally tuned algorithms ready for use 'out of the box'. Our tuning is based on a comprehensive study of parameter variations and options. It would have required several years to complete our study on a single modern desktop computer. Optimal tuning manages to make some very simple and unsexy ideas perform surprisingly well, reducing the need for more ambitious and impressive sounding ones (even if optimally tuned). Our tuning is based on quantitative principles; it can be used for other algorithms as well and implicitly establishes the 'current state of the art' which future proposals may be compared against. It also generates insights previously unavailable about performance comparisons of methods and performance comparisons of different matrix ensembles.

The empirical tuning approach has a larger significance for the field of sparse representations and compressed sensing. Many of the better known papers in this field discuss *what can be proved rigorously*, using mathematical analysis. It requires real mathematical maturity to understand what is being claimed and what the interpretation must be, and to compare claims in competing papers. Often, what can be proved is vague (with unspecified constants) or very weak (unrealistically strong conditions are assumed, far from what can be met in applications). For practical engineering applications it is important to know *what really happens* rather than what can be proved. Empirical studies provide a direct method to give engineers useful guidelines about what really does happen.

The empirical tuning approach also addresses a difficulty many potential users face in addressing the large and growing literature in sparse representations and compressed sensing. In that literature, a rich variety of brand names is developing, where small variations in some already well-known algorithmic schema may lead to the introduction of extravagant new acronyms and phrases. The outsider will be wary of investing the time to digest all this literature and form an accurate understanding of the differences. Empirical tuning efforts group together several differently-named ideas within one family of algorithmic schemas and optimize settings across the whole family, thereby simplifying the situation for many potential users, since the recommended algorithm both runs 'out of the box' and has been tuned to supersede several different earlier proposals.

## II. ITERATIVE ALGORITHMS

Our problem setting will be described with the following notation. An unknown vector $x_0 \in \mathbb{R}^N$ is of interest; we have measurements $y = Ax_0$. Here $A$ is an $n \times N$ matrix and $N > n$. Although the system is underdetermined, it has been shown that, when it exists, sufficient *sparsity* of $x_0$ may allow unique identification of $x_0$. We say that $x_0$ is $k$-sparse if it has at most $k$ nonzeros. In many cases one can exactly recover such a sparse solution $x_0$ as the solution to

$$(P_1) \min \|x\|_1 \text{ subject to } y = Ax,$$

where $\|x\|_1$ denotes the $\ell_1$ norm. This amounts to a large-scale linear programming problem. Unfortunately in some interesting potential applications [31], [32], the matrix $A$ and vector $x_0$ may contain millions of entries and standard linear programming codes may be too slow in those applications. Hence there is widespread interest in finding fast algorithms that work essentially as well; in particular application work by Starck and co-authors [13], [14], [33], [30] and by Elad and co-authors [15], [17] has shown that some very simple iterative algorithms can be strikingly successful on very large problems. In this note we consider two families of such iterative algorithms.

### A. Simple Iterative Algorithms

The first family is inspired by the classical relaxation method for approximate solution of large linear systems. In classical relaxation, one iteratively applies $A$ and its transpose $A'$ to appropriate vectors and under appropriate conditions, the correct solution is obtained as a limit of the process. While the classical theory is inapplicable to underdetermined systems, it has been found both empirically and in theory that a sparsity-promoting variant of relaxation can correctly solve such systems, when they have sufficiently sparse solutions. Starting from $x_1 = 0$, one repeatedly applies this rule:

$$x_{i+1} = \eta_{t_i}(x_i + \kappa \cdot (A' r_i)); \qquad r_i = y - Ax_i;$$

Here $\kappa$ is a relaxation parameter ($0 < \kappa < 1$) and we assume throughout that $A$ is normalized so that its columns have unit length. $\eta_t(\cdot)$ denotes a scalar nonlinearity, applied entrywise; we consider both *hard* thresholding – $\eta_t^H(y) = y 1_{\{|y|>t\}}$ and *soft* thresholding $\eta_t^S(y) = \text{sgn}(y)(|y| - t)_+$. In the above functions $t$ is called the threshold value. Note that if we set $\eta(y) = y$ we would just have classical relaxation. Iterative Soft Thresholding (IST) with a fixed threshold has been used in various settings more than a decade ago – see for example published work of Sylvain Sardy and co-authors [11]. A formal convergence analysis was given by [12] in the determined case. Iterative Hard Thresholding (IHT) was reported useful for several underdetermined cases by Starck, Elad, and their co-authors in papers appearing as early as 2004, [13], [14], [15], [17], [33] often outperforming IST. Other recent examples of such ideas include [16], [19], [18], [34].

These iterative schemes are easy to implement: they require only two matrix-vector products per iteration and some vector additions and subtractions. For certain very large matrices we can rapidly apply $A$ and $A'$ without representing $A$ as a full matrix – examples include partial Fourier and Hadamard transforms. In such settings, the work required scales very favorably with $N$ (eg $N \log(N)$ flops rather than $O(N^2)$).

Actually using such a scheme in practice requires choosing a parameter vector $\theta = (\text{type}, \kappa, t)$; here type = $S$ or $H$ depending as soft or hard thresholding is required; the other parameters are as earlier. Moreover the threshold value $t$ needs to vary from iteration to iteration. The general form in which such schemes are often discussed does not give a true ready-to-run algorithm. This is akin to presenting a cooking recipe listing ingredients for a dish, without listing the needed amounts; it keeps potential users from successfully exploiting the idea.

### B. Composite Iterative Algorithms

In solving determined linear systems, relaxation can often be outperformed by other methods. Because of the similarity of relaxation to IST/IHT schemes, parallel improvements seem worth pursuing in the sparsity setting. A more sophisticated scheme – Two Stage Thresholding (TST) – uses exact solution of small linear systems combined with thresholding before and after this solution. In stage one, we screen for 'significant' nonzeros just as in IST and IHT:

$$v_i = \eta_{t_{1i}}^1(x_i + \kappa A' r_i); \qquad r_i = y - Ax_i;$$

We let $I_i$ denote the combined support of $v_i$ and $x_i$ and we solve

$$w_i = (A'_{I_i} A_{I_i})^{-1} A'_{I_i} y.$$

We then threshold a second time,

$$x_{i+1} = \eta_{t_{2i}}^2(w_i),$$

producing a sparse vector. Here the threshold might be chosen differently in stages 1 and 2 and might depend on the iteration and on measured signal properties. CoSaMP [20] and subspace pursuit [22] may be considered as special cases of TST. This will be clearer when we explain the threshold choice in the next section.

It seems that the use of explicit solutions to the smaller systems might yield improved performance, although at the cost of potentially much more expense per iteration. An important problem is user reticence. In the case of TST there are even more choices to be made than with IST/IHT. This scheme again presents the 'recipe ingredients without recipe amounts' obstacle: users may be turned off by the requirement to specify many such tunable parameters.

### C. Threshold Choice

Effective choice of thresholds is a heavily-developed topic in statistical signal processing. We have focused on two families of tunable alternatives.

*Interference heuristic.* We pretend that the marginal histogram of $A'r$ at sites in the coefficient vector where $x_0(i) = 0$

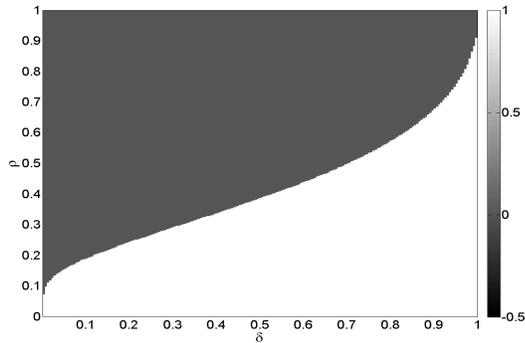
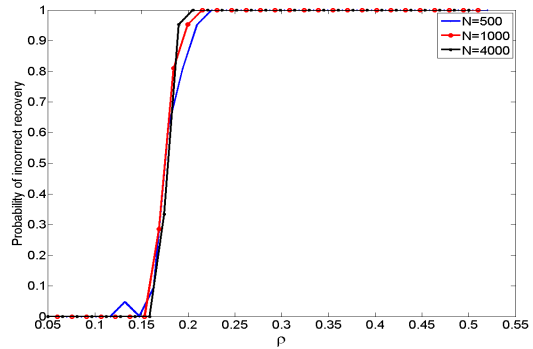

Fig. 1. The success probability of $\ell_1$ in the $(\rho, \delta)$ plane as $N \to \infty$.

Fig. 2. Fraction of unsuccessful recovery attempts by IHT. Here $\delta = .5$ and $\rho = k/n$ is varying. FAR parameter $= 10^{-3}$. Relaxation parameter = 1. At each unique parameter combination, twenty random problem instances were tried. Results are shown at 3 values of $N$: 500, 1000, 4000

is Gaussian, with common standard deviation $\sigma$. We robustly estimate the marginal standard deviation $\sigma$ of the entries in $A'r$ at a given iteration and set the threshold $t$ as a fixed multiple of that standard deviation – $t = \lambda \cdot \sigma$, where $\lambda$ is our chosen threshold control parameter, typically in the range $2 < \lambda < 4$. The underlying rationale for this approach is explained in [10] where its correctness was heavily tested. Under this heuristic, we control the threshold $\lambda$ as in standard detection theory using the False Alarm Rate (FAR); thus $FAR = 2 \cdot \Phi(-\lambda)$ where $\Phi$ denotes the standard normal distribution function.

*Oracle heuristic.* In the TST scheme, imagine that an oracle tells us the true underlying sparsity level $k$, and we scale the threshold adaptively at each iteration so that at stage 1 we yield $\alpha \cdot k$ nonzeros and at stage two $\beta \cdot k$ nonzeros. The method CoSaMP [20] corresponds to $\beta = 2$, $\alpha = 1$, while subspace pursuit [22] corresponds to $\beta = \alpha = 1$.

A problem with the oracle heuristic is that, in interesting applications, there is no such oracle, meaning that we wouldn't in practice ever know what $k$ to use. A problem with the interference heuristic is that the Gaussian model may not work when the matrix is not really 'random'.

## III. PHASE TRANSITIONS

In the case of $\ell_1$ minimization with $A$ a random matrix, there is a well-defined 'breakdown point': $\ell_1$ can successfully recover the sparsest solution provided $k$ is smaller than a certain definite fraction of $n$.

Let $\delta = n/N$ be a normalized measure of problem indeterminacy and let $\rho = k/n$ be a normalized measure of the sparsity. We get a two-dimensional phase space $(\delta, \rho) \in [0, 1]^2$ describing the difficulty of a problem instance – problems are intrinsically harder as one moves up and to the left. Displays indicating success and failure of $\ell_1$ minimization as a function of position in phase space often have an interesting two-phase structure (as shown in figure 1), with phases separated by the curve $(\delta, \rho_{\ell_1}(\delta))$, for a specific function $\rho_{\ell_1}$.

Let $A$ be a random matrix with iid Gaussian entries and let $y = Ax_0$ with $x_0$ $k$-sparse. In [35], [36] one can find explicit formulas for a function $\rho$ definable with the aid of polytope theory and having the following property. *Fix*

$\epsilon > 0$. *The probability that $(P_1)$ recovers the sparsest solution to $y = Ax$ tends to $0$ or $1$ with increasing system size according as $k = n \cdot (\rho_{\ell_1}(n/N) \pm \epsilon)$.* Informally, all that matters is whether $(n/N, k/n)$ lies above or below the curve $(\delta, \rho_{\ell_1}(\delta))$. This is the conclusion of a rigorously proven theorem that describes asymptotic properties as $N \to \infty$; it also describes what actually happens at finite problem sizes [37]. The empirically observed fraction of successful recoveries decays from one to zero as the problem sparsity $\rho = k/n$ varies from just below the critical level $\rho_{\ell_1}(\delta)$ specified in theory to just above it. This transition zone is observed to get increasingly narrow as $N$ increases, matching the theorem, which says that in the large $N$ limit, the zone has vanishing width.

Such sharp phase transitions have also been rigorously proven [38], [39] or empirically observed [10], [6], [40] for other algorithms and/or problem suites. Figure 2 displays behavior of IHT with FAR threshold selection, at a single fixed choice $n/N$ with varying underlying number $k$ of nonzeros. Below a certain threshold, the algorithm works well and above that threshold it fails; the transition zone is narrow, and gets better defined at large problem sizes $N$.

Incidentally, some readers may be unfamiliar with the notion of phase transitions because popular theoretical tools such as coherence [2], [3] and Restricted Isometry Property (RIP) [41] do not really give information about them. It has been shown by Tanner and co-authors [42] that bounds derived from RIP ensure the existence of a region with high success probability in the $\delta$-$\rho$ phase space; however, the actual region is much larger than what those bounds provide.

## IV. ESTIMATING THE EMPIRICAL PHASE TRANSITION

For a given algorithm with a fully specified parameter vector $\theta$, we conduct one phase transition measurement experiment as follows. We fix a *problem suite*, i.e. a matrix ensemble and a coefficient distribution for generating problem instances $(A, x_0)$. For a fixed $N = 800$, we varied $n$ and $k$ through a grid of 900 $\delta$ and $\rho$ combinations, with $\delta$ varying from .05 to 1 in 30 steps and $\rho$ varying from .05 up to a ceiling value $\rho_{max} < 1$ in as many as 30 steps. We then have a grid of $\delta, \rho$ values in parameter space $[0,1]^2$. At each $(\delta, \rho)$ combination, we will take $M$ problem instances; in our case $M = 100$. We also fix a measure of success; see below.

Once we specify the problem size $N$, the experiment is now fully specified; we set $n = \lceil \delta N \rceil$ and $k = \lceil \rho n \rceil$, and generate $M$ problem instances, and obtain $M$ algorithm outputs $\hat{x}_i$, and $M$ success indicators $S_i$, $i = 1, \ldots M$. We run algorithms IHT, IST and TST for 300 iterations.

A problem instance $(y, A, x_0)$ consists of $n \times N$ matrix $A$ from the given matrix ensemble and a $k$-sparse vector $x_0$ from the given coefficient ensemble. Then $y = Ax_0$. The algorithm is called with problem instance $(y, A)$ and it produces a result $\hat{x}$. We declare success if

$$\frac{\|x_0 - \hat{x}\|_2}{\|x_0\|_2} \le \texttt{tol},$$

where `tol` is a given parameter, in our case $10^{-2}$; the variable $S_i$ indicates success on the $i^{\text{th}}$ Monte Carlo realization. We summarize the $M$ Monte Carlo repetitions with $S = \sum_i S_i$. The result $S$ at one experiment will be distributed binomial $\text{Bin}(\pi, M)$ with the success probability $\pi \in [0, 1]$; This probability depends on $k, n, N$ so we write $\pi = \pi(\rho|\delta, N)$.

We measure the location of the phase transition using logistic regression similarly to [6], [43]. The *finite-N phase transition* is the value of $\rho$ at which success probability crosses 50%:

$$\pi(\rho|\delta; N) = \frac{1}{2} \quad \text{at} \quad \rho = \rho^*(\delta; \theta).$$

This notion is well-known in biometrics where the 50% point of the dose-response is called the LD50. (Actually there is a dependence on the tolerance `tol` so $\rho^*(\delta; \theta) \equiv \rho^*(\delta; \theta|N, \texttt{tol})$; but this dependence is found to be weak.

To estimate the phase transition from data, we collect triples $(k, M, S(k, n, N))$ all at one value of $(n, N)$, and model $S(k, n, N) \sim \text{Bin}(\pi_k; M)$ using a generalized linear model with logistic link

$$\text{logit}(\pi) = a + b\rho,$$

where $\rho = k/n$ and $\text{logit}(\pi) = \log(\frac{\pi}{1-\pi})$; in biometric language, we assume that the dose-response probability follows a logistic curve.

The fitted parameters $\hat{a}, \hat{b}$, give the estimated phase transition from

$$\hat{\rho}^*(\delta; \theta) = -\hat{a}/\hat{b}.$$

We denote this estimated value by $\rho^*(\delta; \theta)$ in the rest of the paper.

## V. Tuning Procedure

We conducted extensive computational experiments to evaluate the phase transitions of various algorithms. In all, we performed more than $90,000,000$ reconstructions, using 38 servers at a commercial dedicated server facility for one month. These calculations would have run more than 3 years on a single desktop computer.

For a fixed iterative scheme and a fixed tuning parameter $\theta$, we considered in turn each of several problem suites $\mathcal{S} = (E, C)$, i.e. several random matrix ensembles $E$ and several coefficient amplitude distributions $C$. At each combination, we measured the phase transitions.

In the tuning stage of our project we worked only with the *standard suite* $\mathcal{S}_0$ formed with the Uniform Spherical Ensemble (USE) [1] matrix and constant amplitude distribution on the nonzeros. In the later evaluation stage, other problem suites were used to test the robustness of the tuning. As it turns out, the standard suite is approximately the least favorable case and consequently our tuning works well at all other problem suites.

For a fixed $N = 800$ we measured the empirical phase transition $\rho^*(\delta; \theta)$ as described above. We denote the optimal parameter choice via

$$\theta^*(\delta) = \arg\max_{\theta} \rho^*(\delta; \theta). \quad (1)$$

## VI. Tuning Results

Figure 3 illustrates tuning results for IST on the standard suite $\mathcal{S}_0$. Here $\theta = (RelaxationParameter, FARParameter)$. Panel (a) shows the different optimized phase transitions available by tuning FAR to depend on $\delta$ while the relaxation parameter is fixed. Panel (b) shows the optimally tuned FAR parameters at each given $\delta$ and choice of relaxation parameter. Figures 4 offers the same information for IHT.

Optimum performance of IST occurs at higher values of the false alarm rate than for IHT. Decreasing the relaxation parameter beyond the range shown here does not improve the results for IST and IHT.

Figure 5 illustrates performance of TST for different values of $\theta = (\alpha, \beta)$. Panel (a) shows the different optimized phase transitions available by tuning $\beta$ at fixed $\alpha = 1$ and Panel (b) shows optimal phase transitions with $\alpha = \beta$ varying. Both displays point to the conclusion that $\alpha = \beta = 1$ dominates other choices. Hence subspace pursuit ($\alpha = 1, \beta = 1$) dominates CoSaMP ($\alpha = 1, \beta = 2$).

## VII. Recommended Choices

We provide three versions of iterative algorithms based on our optimal tuning exercise: recommended-IST, recommended-IHT and recommended-TST. They are implemented in Matlab and published at URL `http://sparselab.stanford.edu/ReadyToRun/`.

In our recommended versions, there are no free parameters. The user specifies only the matrix $A$ and the left-hand side $y$. In particular the user does not specify the expected sparsity level, which in most applications cannot be considered known.

These recommended algorithms are not the same as previously published algorithms. For example, recommended TST

---

[1] The columns of these matrices are iid samples from the uniform distribution on the unit sphere in $\mathbb{R}^n$.

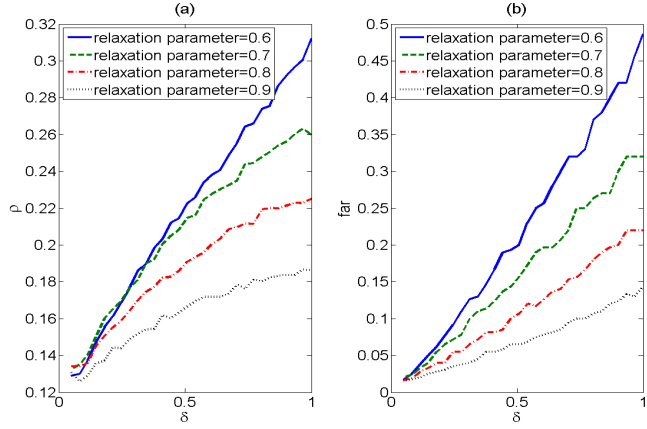

Fig. 3. (a) Optimum phase transitions for IST at several choices of relaxation parameter (b) FAR parameter choice yielding the optimum

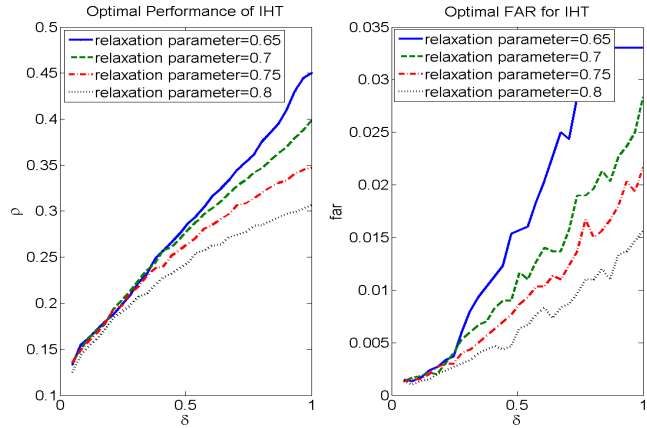

Fig. 4. (a) Optimum phase transitions for IHT at fixed relaxation parameter (b) FAR parameter choice yielding the optimum

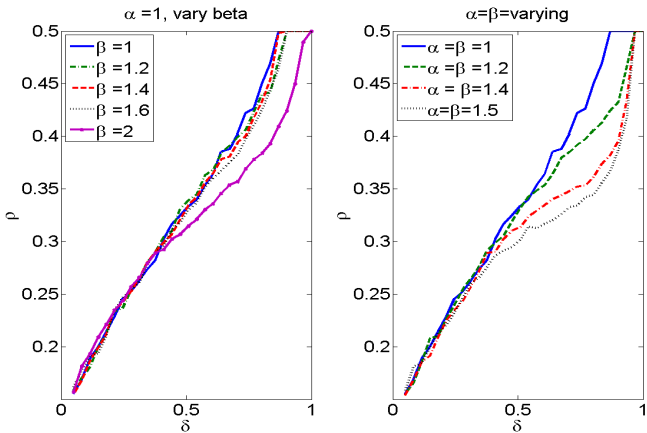

Fig. 5. (a) Empirical phase transitions of TST-$(\alpha,\beta)$ for $\alpha = 1$ and different values of $\beta$; (b) Empirical phase transitions when $\alpha = \beta$.

has parameters $\alpha = 1$ and $\beta = 1$, so it initially seems identical to subspace pursuit [22]. However, subspace pursuit demands an *oracle* to inform the user of the true underlying sparsity of the vector. Recommended-TST has already embedded in it a value for the assumed sparsity level at each $\delta$ (see Table III). If the actual sparsity in $x_0$ is better than the assumed value, the algorithm still works, but if the sparsity is actually worse, the algorithm won't work even if tuned to assume that worse sparsity level. The user does not need to know this number – it is hard-coded. In effect, we have removed the oracle dependence of the subspace pursuit method.

We remind the reader that these algorithms dominate other implementations in the same class. Thus, recommended TST dominates CoSaMP; this is particularly evident for $\delta > .5$ (see Figure 5).

A companion set of algorithms – described later – is available for the case where $A$ is not an explicit matrix but instead a linear operator for which $Av$ and $A'w$ can be computed without storing $A$ as a matrix. Some differences in tuning for that case have been found to be valuable.

We record in the following tables a selection of the optimally tuned parameter values.

TABLE I

RECOMMENDED CHOICES OF FAR AND THE VALUE OF $\rho$ AT WHICH PHASE TRANSITION HAPPENS FOR IST AS A FUNCTION OF $\delta$. OPTIMAL VALUE OF THE RELAXATION PARAMETER $\kappa$ IS $0.6$.

| $\delta$ | .05 | .11 | .21 | .31 | .41 | .5 | .6 | .7 | .8 | .93 |
|---|---|---|---|---|---|---|---|---|---|---|
| $\rho$ | .124 | .13 | .16 | .18 | .2 | .22 | .23 | .25 | .27 | .29 |
| $FAR$ | .02 | .037 | .07 | .12 | .16 | .2 | .25 | .32 | .37 | .42 |

TABLE II

RECOMMENDED CHOICES OF FAR AND THE VALUE OF $\rho$ AT WHICH PHASE TRANSITION HAPPENS FOR IHT AS A FUNCTION OF $\delta$. OPTIMAL VALUE OF THE RELAXATION PARAMETER $\kappa$ IS $0.65$.

| $\delta$ | .05 | .11 | .21 | .41 | .5 | .6 | .7 | .8 | .93 |
|---|---|---|---|---|---|---|---|---|---|
| $\rho$ | .12 | .16 | .18 | .25 | .28 | .31 | .34 | .38 | .41 |
| $100FAR$ | .15 | .2 | .4 | 1.1 | 1.5 | 2 | 2.7 | 3.5 | 4.3 |

TABLE III

RECOMMENDED VALUE OF $\rho$ FOR TST. THE OPTIMAL VALUES OF THE TUNING PARAMETERS ARE $\alpha = \beta = 1$.

| $\delta$ | .05 | .11 | .21 | .31 | .41 | .5 | .6 | .7 | .8 | .93 |
|---|---|---|---|---|---|---|---|---|---|---|
| $\rho$ | .124 | .17 | .22 | .26 | .30 | .33 | .368 | .4 | .44 | .48 |

Figure 6 compares our recommended implementations with each other, with including the theoretical phase transition curve for $\ell_1$ minimization [$\rho_{\ell_1}$; see Section III] and with other algorithms, tuned if appropriate: LARS [4] and OMP [8]. The Figure depicts empirical phase transitions at the Standard Suite. These transitions obey the following ordering:

$$\ell_1 > \text{LARS} > \text{Rec-TST} > \text{Rec-IHT} > \text{Rec-IST};$$

this is exactly the ordering which one would expect based on qualitative grounds; however, it is striking to see how close the curves actually are. On the other hand OMP performance is similar to IHT for $\delta < 0.7$. Both simple iterative algorithms are dramatically less complex to implement and also dramatically cheaper to run on a per iteration basis. It seems that at moderate sparsity levels one would often be satisfied with IHT or IST; particular so for very large problem sizes.

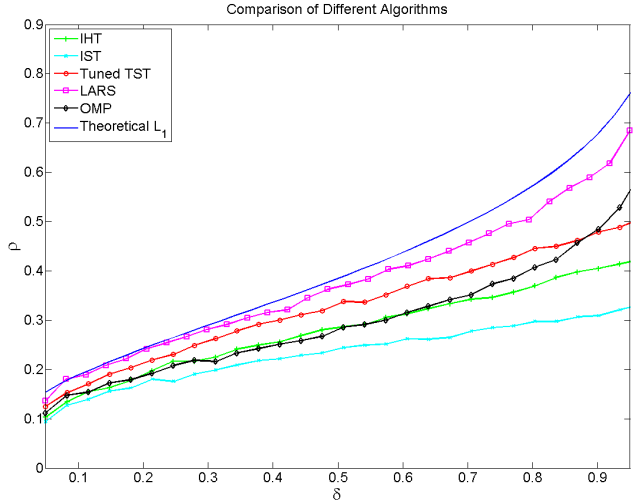

Fig. 6. Phase Transitions of several algorithms at the standard suite. Upper curve: theoretical phase transition, $\ell_1$ minimization; lower curves: Observed transitions of algorithms recommended here.

## VIII. ROBUSTNESS

A *robust* choice of parameters offers a guaranteed level of performance across all situations. Such a choice can be made by solving the maximin problem

$$\theta^r(\delta) = \arg\max_\theta \min_{\mathcal{S}} \rho^*(\delta;\theta;\mathcal{S}).$$

The maximin is achieved at the *least-favorable* suite. As it turns out, our recommended tuning has the maximin property. As described earlier, our tuning results were obtained at the standard suite $\mathcal{S}_0$, with constant amplitude, random-sign coefficients and matrices from USE. We considered a universe of problem suites by varying the matrix ensemble $E$, and the coefficient ensemble $C$. Matrix ensembles included matrices with random $\pm$ entries [Random Sign Ensemble (RSE)] and partial Fourier matrices (to be explained later). We also considered four coefficient ensembles $C$: in addition to the constant amplitude random sign (CARS) ensemble, we considered coefficients from the double exponential distribution, the Cauchy, and the uniform distribution on $[-1, 1]$.

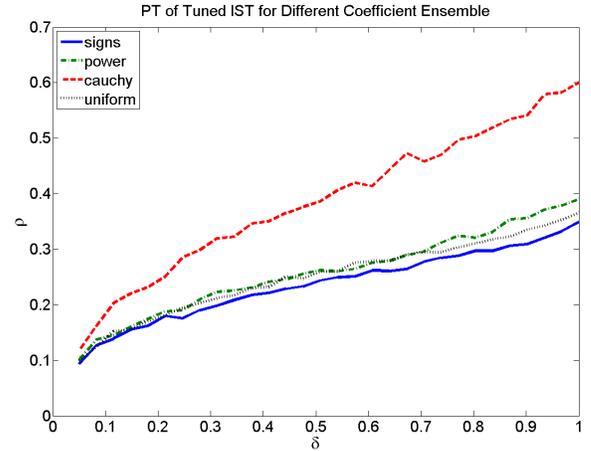

Fig. 7. Observed phase transitions of recommended IST at different coefficient ensembles.

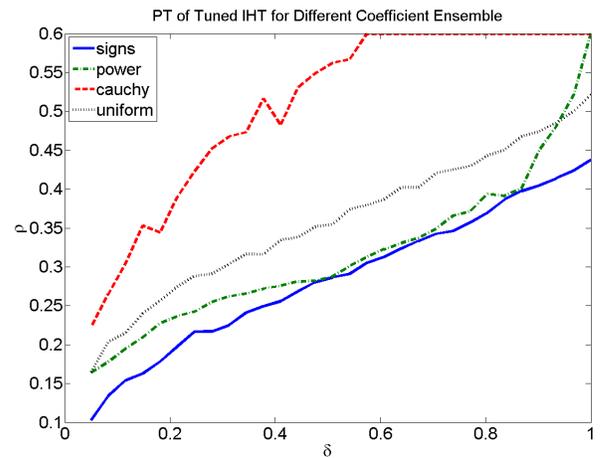

Fig. 8. Observed phase transition of recommended IHT different coefficient ensembles.

Figures 7-8-9 display results for recommended-IST, recommended-IHT, and recommended-TST at a range of problem suites. For all three algorithms, the CARS ensemble is approximately the least favorable coefficient ensemble. Since we have tuned at that ensemble, our choice of tuning parameters can be said to be robust. In other words, if the problem suite is different from the standard suite the phase transition of the algorithm will be above the phase transition we found for $\pm 1$.

Figures 10-11-12 study recommended-IST, recommended-IHT, and recommended-TST at three matrix ensembles USE, Random Sign Ensemble (RSE) where the elements of the matrix are chosen iid from $\pm 1$ and Uniform Random Projection (URP). Results are similar for the RSE and USE ensembles and usually better for URP. A surprising exception to the above pattern is described below in Section X.

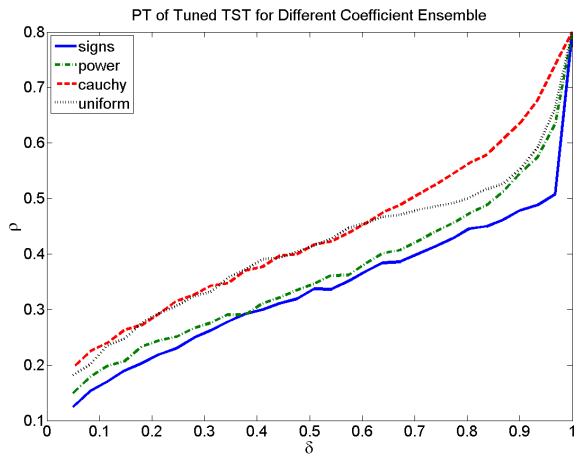

Fig. 9. Observed phase transition of recommended TST for different coefficient ensembles.

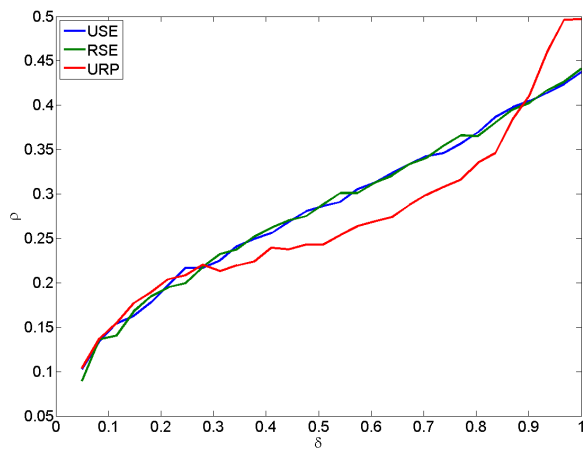

Fig. 11. Observed phase transition of recommended IHT for different matrix ensembles.

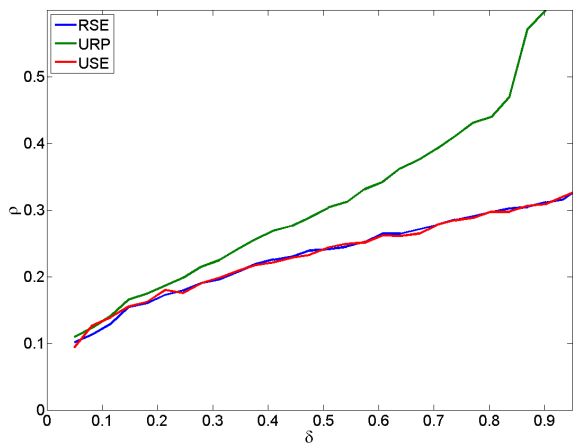

Fig. 10. Observed phase transition of recommended IST for different matrix ensembles.

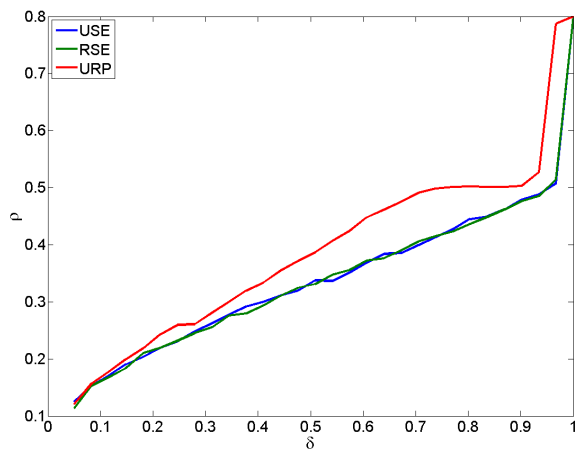

Fig. 12. Observed phase transition of recommended TST for different matrix ensembles.

## IX. Running Times

Algorithm running times were measured on an "Intel 2 Core Processor 2.13GHz" with 3GBytes RAM. All implementations are in Matlab, and in each case the iterative algorithms were run until $\frac{\|y - A\hat{x}\|_2}{\|y\|_2} = .001$.

## X. Ensembles based on Fast Operators

The matrix ensembles discussed so far all used dense matrices with random elements. However, many applications of sparsity-seeking decompositions use linear operators which are never stored as matrices. Typically such operators can be applied rapidly so we call the resulting measurement operators FastOP ensembles. The partial Fourier ensemble [32] provides an example. Here the $n \times N$ matrix $A$ has for its rows a random subset of the $N$ rows in the standard Fourier transform matrix. $Av$ and $A'w$ can both be computed in order $N \log(N)$ time; the comparable dense matrix vector products would cost order $N^2$ flops.

The simple iterative algorithms IHT and IST are particularly suited for use in the FastOp ensemble setting, since they require only repetitive application of $Av$ and $A'w$ interleaved with thresholding.

We considered 1D partial Fourier ensemble, and 1D partial Hadamard ensemble. Figure 13 compares the optimal performance of IHT, IST and TST for partial Fourier matrix ensemble and $\pm 1$ coefficient ensemble.

We found that

- it is beneficial to tune IHT and IST specially for FastOp ensembles, because the previous tuning (aka maximin tuning) was driven by least favorable cases occurring at non-FastOps ensembles. Here such cases are ruled out,

### TABLE IV
AVERAGE RUNNING TIME (SEC) UNTIL $\frac{\|y-A\hat{x}\|_2}{\|y\|_2} = .001$. STANDARD SUITE IS USED IN THESE SIMULATIONS AND THE TIMINGS ARE THE AVERAGE OF 10 INDEPENDENT RUNS.

| $N$ | $\delta$ | $\rho$ | IHT | TST | OMP | LARS |
|---|---|---|---|---|---|---|
| 2000 | 0.9 | 0.17 | 10.6 | 12 | 20 | 28 |
| 4000 | 0.9 | 0.17 | 44.8 | 91.2 | 157 | 216 |
| 6000 | 0.9 | 0.17 | 90 | 286 | 537 | 798 |
| 2000 | 0.7 | 0.28 | 7.2 | 3.3 | 7.6 | 11.5 |
| 4000 | 0.7 | 0.28 | 28.4 | 24.5 | 57.8 | 98.4 |
| 6000 | 0.7 | 0.28 | 64.5 | 118 | 188 | 987 |
| 2000 | 0.5 | 0.2 | 5.8 | 0.91 | 1.5 | 2.7 |
| 4000 | 0.5 | 0.2 | 23 | 7 | 12 | 20 |
| 6000 | 0.5 | 0.2 | 52 | 23 | 38 | 65 |
| 8000 | 0.5 | 0.2 | 91 | 52 | 97 | 164 |
| 10000 | 0.5 | 0.2 | 130 | 100 | 168 | 270 |
| 2000 | 0.3 | 0.12 | 2 | 0.08 | 0.25 | 0.4 |
| 4000 | 0.3 | 0.12 | 9 | 0.65 | 1.8 | 2.6 |
| 8000 | 0.3 | 0.12 | 34 | 5 | 15 | 22 |
| 10000 | 0.3 | 0.12 | 54 | 13 | 28.5 | 38.5 |

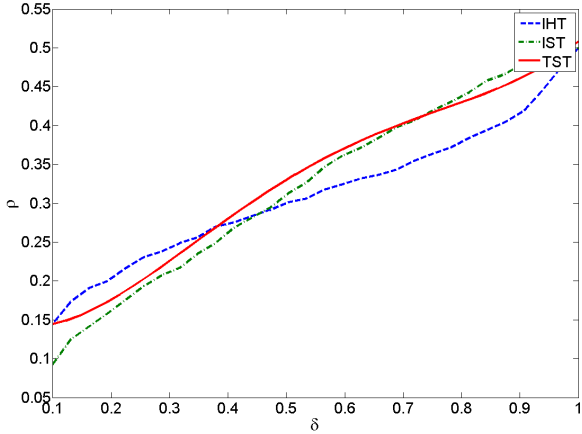

Fig. 13. Comparison of the performance of recommended IHT, IST and TST for partial fourier ensemble.

and maximin tuning only considers a narrower range of relevant cases; the achieved maximin phase transition improves.
- For TST, $\alpha = \beta = 1$ is still optimal, but for the maximin tuning restricted to Fastops $\rho^*$ turns out larger.
- the relaxation parameter in IHT/IST makes essentially no contribution to performance in this setting.
- 1D partial Hadamard and 1D partial Fourier gave very similar results.
- the performance of IHT is very much in line with earlier results for the random matrix ensembles.
- IST behaves dramatically better at partial Fourier ensembles than for the random matrix ensembles (Figure 14) and even outperforms IHT for $\delta > .5$ (Figure 13).

Recommended parameters are shown in Tables V-VI. Running time of the algorithms for partial Fourier are studied in Table VII. As may be noted the execution times of both the fast IHT and fast TST scale favorably with the problem size $N$. In most of our studies TST is faster than IHT and they are both much faster than LARS. The favorable timing results of TST on large problem sizes surprised us.

### TABLE V
RECOMMENDED CHOICES OF FAR AND THE VALUE OF $\rho$ WHERE PHASE TRANSITION HAPPENS FOR IST ALGORITHM. OPTIMAL VALUE OF THE RELAXATION PARAMETER $\kappa = 1$; $E =$ PARTIAL FOURIER

| $\delta$ | .11 | .21 | .31 | .41 | .5 | .6 | .7 | .8 | .9 |
|---|---|---|---|---|---|---|---|---|---|
| $\rho$ | .092 | .16 | .21 | .26 | .31 | .37 | .41 | .44 | .48 |
| FAR | .0209 | .0736 | .13 | .19 | .26 | .32 | .32 | .32 | .32 |

### TABLE VI
RECOMMENDED CHOICES OF FAR AND THE VALUE OF $\rho$ WHERE PHASE TRANSITION HAPPENS FOR IHT ALGORITHM. OPTIMAL VALUE OF THE RELAXATION PARAMETER $\kappa = 1$; $E =$ PARTIAL FOURIER

| $\delta$ | .05 | .11 | .21 | .31 | .41 | .5 | .6 | .7 | .8 |
|---|---|---|---|---|---|---|---|---|---|
| $\rho$ | .056 | .14 | .2 | .24 | .27 | .3 | .32 | .34 | .38 |
| 1000FAR | .3 | .4 | 1.8 | 2.9 | 3.8 | 5 | 5 | 5 | 5 |

### TABLE VII
AVERAGE RUNNING TIME (SEC) UNTIL $\frac{\|y-A\hat{x}\|_2}{\|y\|_2} = .001$. $E =$ PARTIAL FOURIER AND $C =$ SIGNS ARE USED IN THESE SIMULATIONS AND THE TIMINGS ARE THE AVERAGE OF 10 INDEPENDENT RUNS.

| $N$ | $\delta$ | $\rho$ | IHT | TST | LARS |
|---|---|---|---|---|---|
| 8192 | .1 | .1 | .5 | .1 | .6 |
| 16384 | .1 | .1 | 1.2 | .25 | 2.5 |
| 32768 | .1 | .1 | 2.56 | .48 | 10.8 |
| 65536 | .1 | .1 | 8 | 2.3 | 65 |
| 131072 | .1 | .1 | 18 | 5.6 | > 900 |
| 262144 | .1 | .1 | 39 | 13 | > 900 |
| 524288 | .1 | .1 | 85 | 27 | > 900 |
| 16384 | .3 | .18 | .5 | .4 | 25 |
| 8192 | .3 | .18 | .25 | .21 | 5.2 |
| 8192 | .5 | .21 | .18 | .19 | 13.5 |
| 16384 | .5 | .21 | .38 | .4 | 81 |

## XI. DISCUSSION

### A. Before Using These Results

Before applying the results of this project, please note these reminders from the referees to you, the reader:

- Our software already has embedded within it the appropriate values from the tables presented here, so you may not need to copy information from the tables and apply it. However, if you need to code your own implementation, remember that the parameter $\rho = \rho^*$ in our tables specifies the largest workable $k^*$ via $k^* = \lfloor \rho^* \cdot n \rfloor = \lfloor \rho^* \cdot \delta \cdot N \rfloor$.
- Your $A$ matrix must be normalized so that all columns have unit Euclidean norm; for a badly-scaled matrix the algorithms may diverge rapidly.
- The software assumes the sparsity level $k$ is unknown *a priori*, and uses, for each level of the indeterminacy ratio $\delta = n/N$, the largest workable sparsity level $k^*$. If your application provides an oracle that makes $k$ known in advance. you may wish to customize the code to use this information – but this is not necessary.

### B. The Computational Effort-Phase Transition Tradeoff

This project adopted the goal of squeezing the best phase transition performance out of some simple, computationally feasible algorithms. Staring us in the face is the fact that $\ell_1$ minimization generally offers better performance (higher phase transitions) than any of our tuned algorithms.

The referees would like to remind the reader that there is a tradeoff between computational effort and phase transition. At one extreme, explicit combinatorial enumeration, though extravagantly expensive, will obtain optimal phase transition performance. At the other extreme, we have the iterative heuristic algorithms studied here, which can scale to very large problem sizes, but offer more modest phase transition performance. In between, we have convex optimization solvers, which truly solve $\ell_1$ minimization and therefore achieve the phase transition curve $\rho_{\ell_1}$ which, it turns out, dominates the maximin phase transitions of all the iterative methods discussed here. Such solvers run in (pseudo-)polynomial time, but demand considerably more computational effort than the iterative algorithms discussed here.

The tradeoff of computation for phase transition performance is a developing frontier; in particular, there is rapid progress in improving $\ell_1$ solvers. Seven years ago, when the work reported in [13], [14], [15] was conducted, heuristic iterative methods were really the only way to get a start on realistic sized problems. Apparently the popularity of research in Compressed Sensing has led to vast improvements in the performance of available $\ell_1$ solvers. In this project we did not attempt to document what is possible today using such solvers; an interesting follow-up project might involve the tuning of purported $\ell_1$ solvers like Bregman iteration [24] to see if they really can achieve the $\ell_1$ transition; another interesting follow-up project would be to carefully document CPU times of available $\ell_1$ solvers such as L1LS, GPSR, SPGL1, and SpaRSA. Because we follow the paradigm of reproducible computational research it will be possible for others to easily compare such methods with our recommended methods on the exact same range of problem sizes studied in this project.

### C. Contributions

Our work aims to make several contributions:

*1) Helping Potential Users:* The rapidly growing literature on sparsity-promoting reconstruction methods creates difficulties for potential users who would like to apply that knowledge. For a given problem type, there may be many seemingly relevant papers, each promoting specific techniques that have been labeled with catchy branding. When such papers are based on theoretical analysis, the engineer who is a potential user of such ideas may be overwhelmed by the significant amount of mathematical knowledge required to understand and compare the abstract claims being made in such papers. Even when the papers adopt a more familiar engineering approach, the interest of the paper's authors to emphasize the distinctiveness of their algorithm may lead them to use examples or problem settings that are quite different from those used by other authors, making it quite daunting for a newcomer to compare papers and make decisions about what algorithms might conceivably be useful in the user's setting. One easily imagines that in this situation, some users make attempts to digest some of the literature and get inconsistent or disappointing results, not knowing whether the problem is due to misunderstandings, programming errors, or misapplication of techniques, or true limitations of sparsity-promoting methods.

Our results may help potential users more easily evaluate the potential benefits of well-chosen sparsity-promoting reconstruction methods. Several algorithms are made freely available, their properties are carefully described and compared on a common basis, and the underlying problem suites and performance metrics are available for careful study. We give the user an ability to get started with these ideas very directly and transparently,

*2) Defining the State of the Art:* Our work may be useful to active researchers in the field of sparsity-promoting methods. Our study effectively defines the *current state of the art* (CSA), a precisely-specified set of quantitative performance standards which are the best we currently know how to do (within a certain class of algorithms). Once this is defined, any newly proposed algorithm can be evaluated with reference to the CSA, and other researchers can use this comparison to understand the relative improvement, if any, offered by the proposal. Over time, as genuine improvements emerge, the CSA will evolve, by definition always offering the best known current performance.

More broadly, a researcher with a new method *not* intended for comparison with the CSA on a standard suite may expand the set of suites beyond those studied here, adding to such a suite a new matrix or coefficient ensemble, or may provide a new measure of success. This allows the researcher to clearly demonstrate for colleagues the arena where the method is intended to contribute.

The effort to create standard performance metrics and define CSA performance has been valuable in many fields of image and signal processing. Indeed, one can argue that in fingerprint recognition and face recognition, the moment when those fields really started to make progress is precisely the moment when defined databases and success metrics were made available for community use, allowing systematic comparison of algorithms [44], [45]. Similarly, the regular publication of standardized challenge problems in arenas like protein structure prediction is said to have utterly transformed the field. Defining a CSA for sparsity-promoting methods can likewise be expected to lead to much faster and more reliable progress.

Standardization has another predictable effect: it may *improve communication among researchers in the field* reducing the impact of marketing, branding and prestige and increasing the focus on objective measures performance. The example of protein structure prediction (CASP) bears this out [46]; before CASP, certain approaches to protein structure prediction were considered more likely to work than others, but it has been reported that CASP upset expectations, reversing the order of preference for certain algorithms [47]. In fact our work may already be showing this effect: one of our findings is that the more prominent algorithm CoSaMP is dominated by the less well-known algorithm Subspace Pursuit at the suites we have studied and for the success measure we have used.

*3) Promoting Reproducible Computational Research:* An implicit but still important contribution of our work is adherence to the paradigm of Reproducible Computational Research [48]. The data and code required to reproduce our results are freely available – and not just the conclusions. A researcher or potential user can study our implementation and tuning of an algorithm or our definitions of performance metric, or our collection of problem suites. This speeds up progress in developing and validating new sparsity-seeking algorithms.

- Researchers who feel our study ignores important aspects of performance can easily define new metrics of success which better reflect their views of what is important, and to conduct parallel studies with that new metric.
- Potential users interested in a specific problem suite we haven't studied, but which is of direct interest in their application, may easily extend our software to add a new suite to the available collection and then run a robustness study or even a tuning study focused on that suite.

By sharing the code underlying our study, we promote further developments of new applications and new methods which outperform current ones.

## XII. Conclusions

We defined a set of problem suites and two algorithmic schemes that cover several distinctly branded methods in the literature. We defined algorithm performance using the notion of empirical phase transition and made millions of reconstruction attempts, while systematically varying problem specifications. We identified specific parameter choices which are optimal at so-called standard suites, for specific sparsity-indeterminacy combinations. This produced recommended- IST, IHT and TST algorithms, coded in Matlab and are freely available at URL `sparselab.stanford.edu/ReadyToRun`. They can be used 'out of the box' on problem instances of the type we have studied; the user need not specify any parameters whatever in order to run them; simply providing the matrix $A$ and the left-hand side $y$ of the system $y = Ax$.

Our studies included extensive computations at other suites besides the standard one, verifying the robustness of our parameter choices and checking that at those other suites the recommended algorithms generally behave better than they do at the standard suite.

The standard suite involves random matrices, but many applications of sparsity-seeking algorithms, particularly in compressed sensing, use structured matrices. We did consider an important class of structured matrices based on fast transforms such as 1D and 2D Fourier transforms and fast Hadamard transforms. Such matrices admit of rapid computations with very large problem sizes and in some cases are actually demanded by the application, for example in MR imaging and NMR spectroscopy. We studied some very large problem sizes with problem suites based on these fast transforms and found results largely matching those we found on the random ensembles. In one case – IST with 2D partial fourier – we found tunings which generate unexpectedly high phase transitions markedly better than what we saw for the standard suite. We published recommended choices of IHT, IST for use with such ensembles defined by those fast operators.

We reached the following empirical findings at the 'random' matrix ensembles

- Phase transitions for optimally-tuned algorithms obey the ordering recommended-TST > recommended-IHT > recommended-IST.
- Setting the relaxation parameter to 0.6 for IST and 0.65 for IHT improves performance of those algorithms significantly. Relaxation has no noticeable effect on performance of TST.
- Performance of the matrix ensemble USE (start with iid Gaussian entries then normalize column lengths) is very similar to RSE (random ±1 entries).
- The distribution of coefficient amplitudes in the solution $x_0$ matters very much to these algorithms. The worst case is when all nonzeros have the same amplitude.
- For a given problem suite, the number of iterations to reach a given accuracy of recovery does not seem to depend on problem size, except perhaps near to phase transition.
- Subspace pursuit works better than CoSaMP on the standard suite. Our recommended TST algorithm is essentially subspace pursuit, without need of an oracle.

Our conclusions for the case where the matrix is defined using a fast linear operator were listed in Section X.

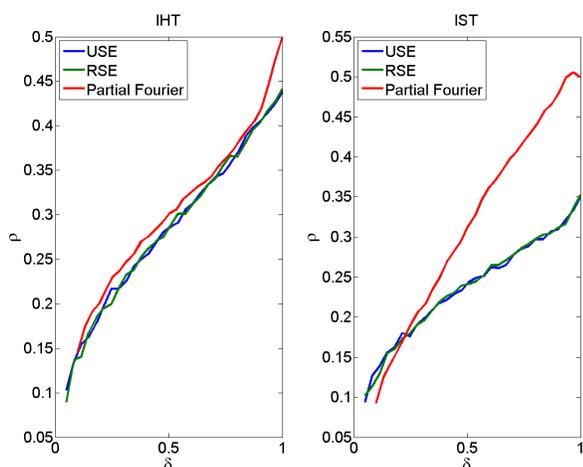

Fig. 14. (a) Phase transitions of recommended IHT for different matrix ensembles (b)Phase transitions of recommended IST for different matrix ensembles

ACKNOWLEDGEMENTS

Thanks to Jared Tanner for valuable suggestions on an early draft, and the Special Issue Editors, especially Rick Chartrand, for their service to the IEEE. This work was partially supported by NSF DMS 05-05303.